\documentclass[a4paper,twosided]{article}

\usepackage[T1]{fontenc}
\usepackage[utf8]{inputenc}
\usepackage{amsmath}
\usepackage{amssymb}
\usepackage{amsthm} 
\usepackage{subfig}
\usepackage{graphicx}
\usepackage[margin=2cm]{geometry}
\usepackage[
  justification=centerlast,
  singlelinecheck=true,
  font={small,it},
  labelfont=bf
]{caption}[2004/07/16]

\newtheorem{theorem}{Theorem}

\newcommand{\dfnd}{\overset{\text{def}}{=}}
\newcommand{\maxtrigdeg}{$2KR_I$}
\newcommand{\RII}{ R_{I\!I} }
\newcommand{\firstProg}{ \texttt{expl\_trasf.c} }
\newcommand{\lieDer}[1]{ \{#1\} }
\newcommand{\delenda}[1]{ F_{#1} }
\newcommand{\hatdel}[1]{ \hat{ \delenda{#1} } }
\newcommand{\ham}[1]{ F_{#1} }
\newcommand{\hatham}[1]{ \hat{ \ham{#1} } }
\newcommand{\hamavg}[2]{\langle F_{#1}^{#2} \rangle}
\newcommand{\pluseq}{ \,\mathtt{+\!=}\, }
\newcommand{\hamcap}{H}
\newcommand{\Pphi}{P}
\newcommand{\ecap}{\mathcal{E}}
\newcommand{\zcap}{\zeta}
\newcommand{\acap}{a}
\newcommand{\hbound}[3]{ \mathcal{F}_{#1}^{(#2,#3)} }
\newcommand{\hathbound}[3]{ \hat{\mathcal{F}}_{#1}^{(#2,#3)} }
\newcommand{\vbound}[3]{ \mathcal{F}_{#1}^{(#2,#3)} }
\newcommand{\hatvbound}[3]{ \hat{\mathcal{F}}_{#1}^{(#2,#3)} }
\newcommand{\avgbound}[2]{ \mathcal{B}^{(#1,#2)} }
\newcommand{\genbound}[3]{ \mathcal{G}_{#1,#2}^{(#3)} }

\newcommand{\setpol}[2]{ \mathbb{P}_{#1,#2} }

\title{Hamiltonian Control of Magnetic 
Field Lines: Computer Assisted Results 
Proving the Existence of KAM Barriers}
\author{Lorenzo Valvo\ \& Ugo Locatelli$\thanks{\texttt{locatell@mat.uniroma2.it}}$ \\
  {\small Dipartimento di Matematica } \\
  {\small Universit\`a degli Studi di Roma ``Tor Vergata''}\\
  {\small Via della Ricerca Scientifica 1} \\
  {\small 00133 Roma, Italy }
 }

\begin{document}

\maketitle

\begin{abstract}
  We reconsider a control theory for Hamiltonian systems, that was
  introduced on the basis of KAM theory and applied to a model of
  magnetic field in previous articles. By a combination of Frequency
  Analysis and of a rigorous (Computer Assisted) KAM algorithm we
  prove that in the phase space of the magnetic field, due to the
  control term, a set of invariant tori appear, and it acts as a
  transport barrier. Our analysis, which is common (but often also
  limited) to Celestial Mechanics, is based on a normal form approach;
  it is also quite general and can be applied to quasi-integrable
  Hamiltonian systems satisfying a few additional mild assumptions. As
  a novelty with respect to the works that in the last two decades
  applied Computer Assisted Proofs into the framework of KAM theory,
  we provide all the codes allowing to produce our results. They are
  collected in a software package that is publicly available from the
  {\it Mendeley Data} repository. All these codes are designed in such
  a way to be easy-to-use, also for what concerns eventual adaptations
  for applications to similar problems.
\end{abstract}

In a few works based on variants 
of our CAP (\cite{giorgilli_kolmogorov_2009},
\cite{locatelli_invariant_2000},
\cite{gabern_construction_2005}) it 
was evident that performing a few 
preliminary manipulations on the 
initial Hamiltonian can really 
improve the final outcome of the 
CAP. In this same spirit we are
introducing the
Hamiltonian $\hamcap^{(0)}$, from
which the whole CAP will be started.
A key topic of research in plasma physics is the 
confinement of the plasma itself~\cite{hazeltine_plasma_2003}. 
One popular technique in laboratory experiments and
industrial applications is 
magnetic confinement: charged particles are
confined by the application of strong magnetic
fields in suitable configurations; the question 
then becomes which are the most effective 
field geometries. A necessary (but unfortunately
not sufficient) condition for the confinement
of the plasma is the confinement of the magnetic 
field itself: in fact it is known that the
trajectories of the charged particles approximately 
follow the magnetic field lines 
\cite{northrop_adiabatic_1963}.

One common magnetic configuration is a 3 dimensional 
magnetic field whose lines are fold on 2 dimensional
toroidal surfaces, symmetric by rotations around a
central axis (\textit{axisymmetry}). The lines of 
this type of field can be seen as the trajectories 
of a one-and-a-half degrees of freedom 
(henceforth, $1\tfrac{1}{2}$ d.o.f.)
Hamiltonian system \cite{abdullaev_construction_2006}. 
The integrability of this Hamiltonian
system corresponds to the confinement of the magnetic 
field. In presence of a perturbation this
property is broken: the magnetic field changes
its topology and we observe the outbreak of
chaos. This type of magnetic perturbations
are named \textit{tearing modes} in plasma physics
(see for instance \cite{bellan_fundamentals_2008},
chapter 12).

One possibility to reduce chaos is control
theory. We are interested here in an 
original approach 
developed by M. Vittot and coauthors in a 
series of papers (we recall in particular 
\cite{vittot_perturbation_2004},
\cite{ciraolo_control_2004},
\cite{vittot_localised_2005}), and applied 
to the problem of magnetic field lines in
\cite{chandre_control_2006}. Assuming to 
have an unpertubed hamiltonian $H_0$ 
and a perturbation $v$, they 
show how to build a control term $f$ of order
$v^2$ such that the Hamiltonian
$H_0 + v + f$ has an invariant torus in 
phase space. In our case the phase space has 
3 dimensions, reduced to 2 by conservation of
energy, so the invariant curve acts as a 
barrier (commoly referred to as 
\textit{transport barrier} in plasma physics)
confining (part of) the trajectories.  

In this work we want to investigate a few 
aspects that were not considered (at least, 
not extensively) in \cite{chandre_control_2006}.
First of all, it is not clear that the 
invariant curve is a torus, and if this 
is the case, which frequency is associated to it.
Second,  we try 
to measure quantitatively the actual improvement 
in the integrability of the system due to the 
control term $f$. Third, the control theory
needs a localization in phase space (as in KAM 
theory); but it is not clear if the invariant 
curve is located at that same point of phase
space (as in KAM theory), nor if it is a
single curve, or if the barrier has an actual
``thickness'' (as in KAM theory). 

To answer these questions we first apply the Frequency Analysis (as
introduced by Laskar, see, e.g.,~\cite{laskar_introduction_1999}) to
deduce some qualitative properties of the system, and how these
properties change as the control is applied.  Then, we also show the
existence of invariant KAM tori through a Computer Assisted Proof
(henceforth CAP). This second result is of quantitative nature because
it employs validated numerics:
we replace any number by an upper bound, 
or a lower bound, or both (in this case 
we speak of \textit{interval arithmetics}, 
because a number is replaced by an interval) 
to take into account any possible error 
in the determination of that number. Also
the operations among numbers are modified 
accordingly. The error can have any origin: 
in this work it coincides with the error 
made by the computing machine in representing 
a real number by a floating point number; 
however, if the number has a physical meaning, 
the error may be the maximum precision allowed 
from the experimental setup.

This paper is organised as follows:
in section \ref{sec:model} we describe the model;
in section \ref{sec:afma} we report the results
of the Frequency Analysis; in 
section \ref{sec:cap} we describe the CAP; 
in section \ref{sec:appo} we apply
the CAP to the problem of
the magnetic field. A final section 
\ref{sec:conclusions} is dedicated to 
the conclusions.

\section{The model}
\label{sec:model}

A three dimensional azimuthally symmetric 
magnetic field can be 
represented as a $1\tfrac{1}{2}$ d.o.f. 
Hamiltonian system
\cite{abdullaev_construction_2006}. The toroidal 
angle $\varphi\in[0,2\pi)$ has the role of time, while
the toroidal flux $\psi\in[0,1]$ and the poloidal
angle $\theta\in[0,2\pi)$ are the conjugated action and 
angle variables, respectively. The poloidal
flux is the Hamiltonian, hence denoted by $H$. 
At equilibrium, the poloidal flux can be 
written as a function of the toroidal flux 
alone and the dynamical system is integrable: 
the magnetic field lines are confined. 
Traditionally, we write
$H(\psi)\,=\,\int 1/q(\psi) d\psi$, where the
function $q$ is called the \textit{safety factor}.
A common model, described in 
\cite{abdullaev_construction_2006} 
and also adopted in 
\cite{chandre_control_2006}, is defined by setting
\begin{equation}
  \frac{1}{q(\psi)}\,=\,\frac{(2-\psi)(2-2\psi+\psi^2)}{4}
\end{equation}
so that 
\begin{equation}
  \label{eq:hamiltonian1}
  H(\psi)\,=\,\psi\,-\tfrac{3}{4}\psi^2\,+\,\tfrac{1}{3}\psi^3\,
    -\,\tfrac{1}{16}\psi^4. 
\end{equation}

The perturbation $v$ introduces a dependence 
of the poloidal flux $H$ on the angles 
$\theta$ and $\varphi$, and the integrability
is lost. As we also recalled in the
introduction, a perturbation which
changes the topology of the magnetic field
lines is called a \textit{tearing mode} in 
plasma physics. As in reference 
\cite{chandre_control_2006}, we set
\begin{equation}
  \label{eq:perturbation0}
  v(\theta,\varphi)\,=\,\varepsilon\big(
    \cos(2\theta-\varphi)\,+\,\cos(3\theta-2\varphi)\, \big).
\end{equation}
The dynamical functions defined on the phase space 
can be grouped into the sets $\setpol{l}{s}$ 
defined as follows
\begin{equation}
  g\in\setpol{l}{s}\,\iff\,
    g\,=\,\sum_{ \max(|k_1|,|k_2|)\le s }
      \,g_{l,k_1,k_2}\,\psi^l\,e^{ik_1\theta+ik_2\varphi}
 \label{eq:gp}
\end{equation}
and the coefficients $g_{l,k_1,k_2}$ are sometimes 
called the \textit{Taylor-Fourier coefficients} of $g$.
If we introduce the constants $l_\text{max}$ and $K$, 
respectively equal to the highest polynomial degree of 
$H(\psi)$ and to the trigonometric degree of the 
perturbation $v$, from equations \eqref{eq:hamiltonian1} 
and \eqref{eq:perturbation0} we have 
\mbox{$H\in\bigcup_{l=0}^{l_\text{max}}\setpol{l}{0}$} 
with $l_\text{max}=4$
and \mbox{$v\in\setpol{0}{K}$} with $K=3$. 

We look for a torus around $\Psi\,=\,0.35$, 
which is located midway between the two resonant
surfaces \mbox{$\psi_{3,2}\approx 0.266$} and 
\mbox{$\psi_{2,1}\approx 0.456$} (implicit solutions of
the equations $q(\psi_{3,2}) = 3/2$ and $q(\psi_{2,1}) = 2$).
So we make a translation of the $\psi$ coordinate,
$\psi\,\mapsto \psi-\Psi$. We also authonomize the
system by introducing a dummy momentum coordinate
$\Pphi$
conjugated to $\varphi$. The perturbed Hamiltonian
$H+v$ becomes
\begin{equation}
 \tilde{H}\,=\,\tilde{\omega}\psi\,+\,\Pphi\,+\,
  h(\theta,\varphi,\psi)\,+\,v(\theta,\varphi)\,
  ,\qquad 
   h(\theta,\varphi,\psi)\,=\,
    \sum_{l= 2}^{l_\text{max}} h_l.
  \label{eq:hamiltonian2} 
\end{equation}
Note that in equation \eqref{eq:hamiltonian2} we
considered a dependence of $h$ on $\theta$ and 
$\varphi$ for the sake of generality, although
in our example this is not the case. 

To build the control term we need to solve
the \textit{homological equation} 
\begin{equation}
  (\tilde{\omega}\partial_\theta + \partial_\varphi)X\,+\,v(\theta,\varphi)\,= 0.
\label{eq:homological1}
\end{equation}
In \cite{chandre_control_2006} the
authors introduce the operator $\Gamma$ which
is the fundamental solution of the homological
equation. So we can write
\begin{equation}
  X\,=\,-\,\Gamma_{\tilde{\omega}}\,v(\theta,\varphi)\, 
   \equiv\,\sum_{ (k_{1}\,,\,k_{2}) \in \mathbb{Z}^2\setminus\{(0\,,\,0)\} }\,
    -\frac{ v_{k_1,k_2} }{k_1\tilde{\omega} + k_2}\, 
      e^{ik_1\theta + i k_2\varphi}.
  \label{eq:defGamma}
\end{equation}
The solvability of the homological equation
is granted by the assumption that $\tilde{\omega}$ 
is \textit{non-resonant}: for any
$(k_{1}\,,\,k_{2}) \in \mathbb{Z}^2\setminus\{(0\,,\,0)\}$
such that \mbox{$0<\max(|k_1|,|k_2|)\leq K$},
we have \mbox{$\big|k_1 \tilde{\omega} + k_2\big|\, \neq\,0$}.

\begin{figure}
  \centering
  \includegraphics[width=0.45\textwidth]{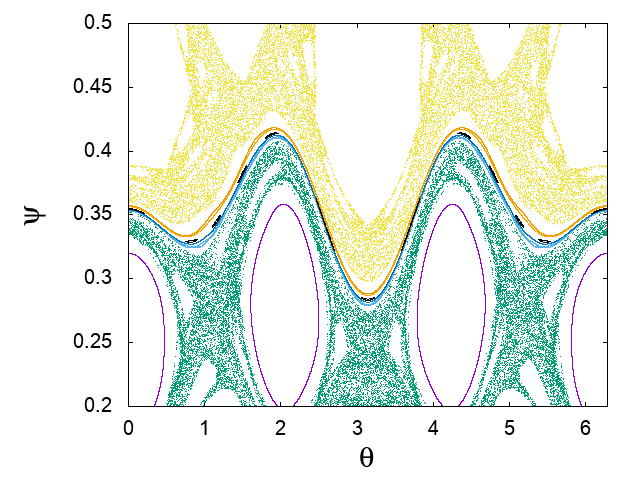}
  \caption{ Phase portraits given by the time-$2\pi$ mappings for the
    \label{fig:phase_portrait}
    Hamiltonian $H+v+f$ (see equations \eqref{eq:hamiltonian1},
    \eqref{eq:perturbation0}, \eqref{eq:control3}). The equations of
    motion were numerically integrated by using a leap-frog method in
    the case with $\varepsilon=0.003$.  }
\end{figure}

The control term for the Hamiltonian 
\eqref{eq:hamiltonian2} is
\begin{equation}
 \label{eq:control2} 
   f\,=\,-\frac{\varepsilon^2}{2}{\lieDer{X}}^2 h_2
    \,=\,-\,\frac{1}{2}\,\frac{d}{d\psi}
     \Big(-\frac{1}{q(\psi)}\Big)\Big|_{\psi=\Psi}\,
      \big(\,\partial_\theta \Gamma_{\tilde{\omega}} v\,\big)^2
\end{equation}
where we have introduced the Lie bracket
\begin{equation}
  \lieDer{F}\,G\,\dfnd\,\partial_\theta F\,\partial_\psi G\,
   +\,\partial_\varphi F\,\partial_{\Pphi} G\,-\,
    \partial_\theta G\,\partial_\psi F\,-\,
     \partial_\varphi G\,\partial_{\Pphi} F
\label{eq:lieBracket}
\end{equation} 
among any two dynamical functions $F$ and $G$.
If we plug equations \eqref{eq:hamiltonian2}, 
\eqref{eq:perturbation0} and \eqref{eq:defGamma} 
into \eqref{eq:control2}, we get 
\begin{equation}
    f\,
       =\varepsilon^2\,\big[ -\tfrac{9}{2}\cos(6\theta - 4\varphi)\, +\,
         6\,\big( \cos(\theta-\varphi) - \cos( 5\theta-3\varphi) \big)  
         +\,2\,\cos( 4\theta-2\varphi)  \big] \partial_{\psi\psi}^2 h_2.
  \label{eq:control3}
\end{equation}
Let us also recall that the third expression 
of equation \eqref{eq:control2} coincides 
with equation (14) of \cite{chandre_control_2006}.

The plots of some trajectories of the dynamical system with
Hamiltonian $H+v+f$ in Figure~\ref{fig:phase_portrait}.  It is
possible to appreciate the existence of a ``thin'' set of trajectories
with initial conditions around the point
$(\psi\,,\,\theta)=(0.35\,,\,0)$ that look to lie on invariant curves.
We stress that the Hamiltonian $H+v$ is not localized around this value,
but the control term $f$ was computed around it.

\section{Qualitative Analysis}
\label{sec:afma}

In this section we apply the Frequency
Analysis of Laskar 
\cite{laskar_introduction_1999} to 
study how the control term affects the 
dynamics induced by the Hamiltonian $H$
of equation \eqref{eq:hamiltonian1}.
In particular we can characterize a 
Hamiltonian system by associating to an 
initial datum for the action (here $\psi_0$)
the frequency $\omega$ of the corresponding 
trajectory; by repeating for many values 
of $\psi_0$ we build a graphics of 
$\omega \text{ vs } \psi_0$, called the 
Frequency Action Map (henceforth FAM). 
To perform the association, we solve numerically
the dynamics for a given $\psi_0$ and we determine 
$\omega$ by applying the fundamental integral 
of the frequency analysis 
\cite{laskar_introduction_1999}. This is 
repeated many times. 

\begin{figure}
\centering
\subfloat[][
  \label{fig:afm1}
]{
  \includegraphics[width=0.45\textwidth]{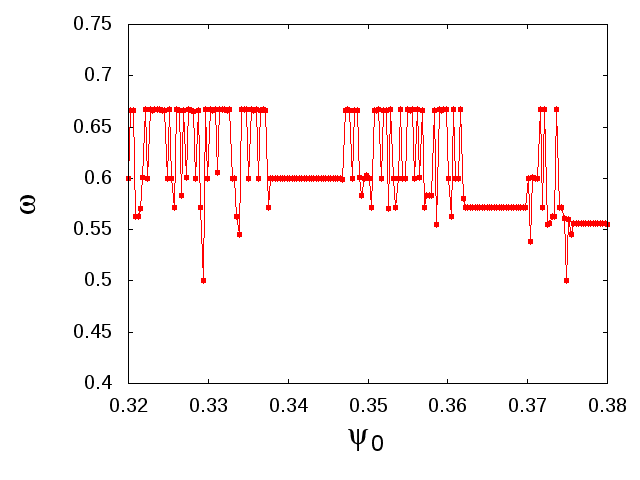}
} 
\hspace{0.05\textwidth} 
\subfloat[][
  \label{fig:afm2}
]{
  \includegraphics[width=0.45\textwidth]{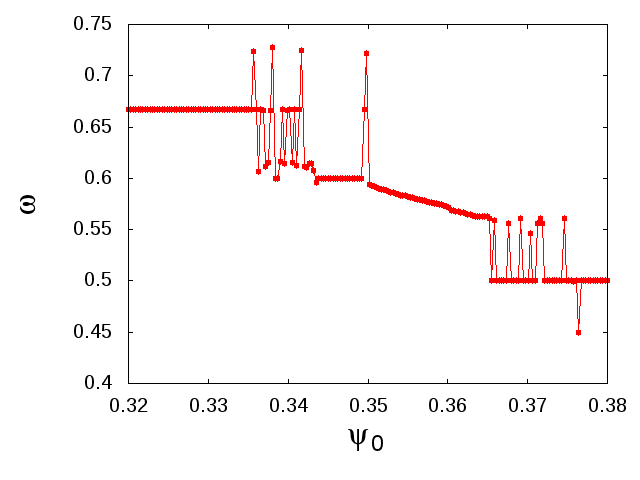}
} 
\caption{
  \label{fig:afm}
  Two FAMs for the Hamiltonian 
  system corresponding to the Hamiltonians
  $H+v$ (on the left) and $H+v+f$  
  (see equations \eqref{eq:hamiltonian1},
  \eqref{eq:perturbation0} and
  \eqref{eq:control3}), for 
  $\varepsilon=0.0012$. 
  In both cases we chose 200 equidistributed
  values of $\psi_0$, and for each of them 
  we run a simulation of 
  \mbox{$(2^{15}+1)$} perturbation periods.
  The initial values of the 
  variables $\theta$ and $\varphi$ were always 
  set to $0$. The equations of motion were 
  solved by a symmetric splitting
  method of order two. 
}
\end{figure}

First we consider the dynamics determined by
the Hamiltonian $H+v$ with $H$ and $v$ as in 
equations \eqref{eq:hamiltonian1} and
\eqref{eq:perturbation0}. We set 
$\varepsilon=0.0012$ and we compute 
the map relative to the interval
$0.32<\psi_0<0.38$, evenly distributed around 
$\Psi=0.35$. Then we add the control term
$f$ given by equation~\eqref{eq:control3} and
we repeat the analysis. The two FAMs are
reported in figures~\ref{fig:afm1}
and~\ref{fig:afm2}. In the
case without control term the graphics is dominated by scattered
points, which correspond to chaotic regions, and flat branches,
corresponding to resonances.  On the right (in presence of the control
term), the graphics is also dominated mostly by resonances, but there
is also a regularly decreasing region (where the behavior looks
linear), approximately for $0.35\leq \psi_0\leq 0.365$. According to
the basics of the FAM method, this corresponds to a zone of the phase
space that is filled with invariant (KAM) tori.  We observe that, even
though we localized the Hamiltonian around $\Psi=0.35$ to compute the
control term, the invariant tori do not appear around $\Psi$, but they
are displaced to its right. This is usually seen also in KAM theory.

\begin{figure}
\centering
\includegraphics[width=\textwidth]{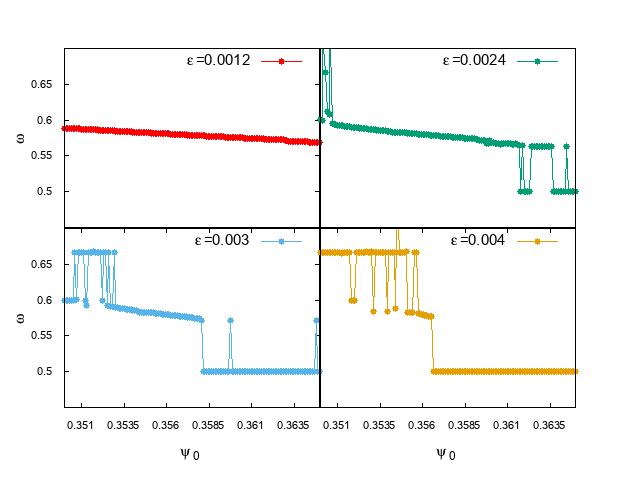}
\caption{
  \label{fig:confrontation}
  The action-frequency map for the system $H+v+f$
  at different values of $\varepsilon$, in the 
  region where the KAM tori appeared as a 
  consequence of the control term. For each
  value of $\varepsilon$ we numerically 
  computed $150$ trajectories of 
  $(2^{15}+1)$ perturbation periods.
}
\end{figure}

Next we computed different FAMs on the 
interval $0.35\leq \psi_0\leq 0.365$,
increasing the value of $\varepsilon$.
The results are reported in the four
boxes of figure~\ref{fig:confrontation}.
As $\varepsilon$ is doubled from 
$0.0012$ to $0.0024$, the width of the 
regular pattern is left nearly unchanged; 
then it slowly shrinks, as $\varepsilon$ 
increases, leaving space for more 
chaos and resonances. At 
$\varepsilon=0.003$ the width of the
barrier is more than halved, while at 
$\varepsilon=0.004$ (the value 
considered in \cite{chandre_control_2006}) 
it is extremely small.
Two similar plots (that we omit for the sake
of brevity) allows us to fix to
$\varepsilon^*\in(0.0045\,,\,0.0046)$ the
threshold value for the disappearance of
the invariant tori separating the chaotic
regions. In particular, for
$\varepsilon=0.0046$ the results provided
by the FAM (in the case of the controlled
Hamiltonian $H+v+f$) have the same look
with respect to those reported in
figure~\ref{fig:afm1} (corresponding to
the uncontrolled case $H+v$).

\begin{figure}
\centering
\includegraphics[width=0.7\textwidth]{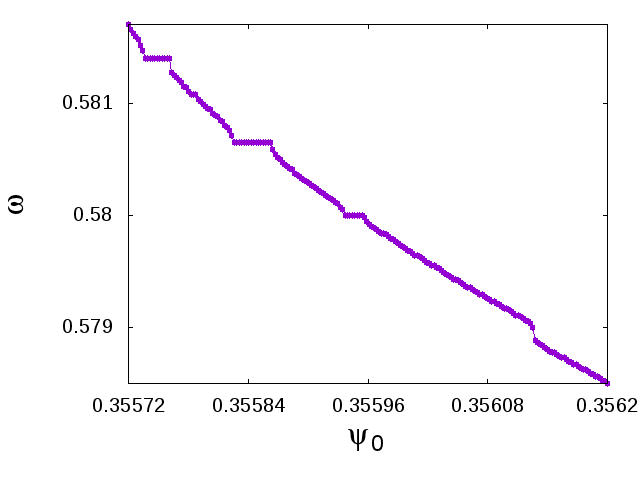}
\caption{
  \label{fig:afm3}
  Action-frequency map analysis for a region
  of phase space of the system $H+f$ filled
  with invariant tori, for $\varepsilon= 0.004$.
  The graphics is made of $200$ points.
  For each point we run a simulation of 
  \mbox{$(2^{16}+1)$} perturbation periods.
}
\end{figure}

Finally, we perform a further magnification
of the interval 
$0.35572 \leq \psi_0 \leq 0.3562$ in the
case $\varepsilon=0.004$; this is 
reported in figure~\ref{fig:afm3}. At
this value of the perturbation 
parameter we have a good compromise 
among being close to the breaking
threshold of the invariant tori, and 
being able to distinguish the linear 
region of the FAM. 
In figure~\ref{fig:afm3} there are 
four small plateau, corresponding 
to the (resonant) frequencies   
\begin{equation}
  \omega_1\,=\,\frac{25}{43}\,\approx\,0.5813968\,
  ,\quad\omega_2\,=\,\frac{43}{74}\,\approx\,0.5810763\,
  ,\quad\omega_3\,=\,\frac{18}{31}\,\approx\,0.5806491\,
  ,\quad\omega_4\,=\,\frac{29}{50}\,\approx\,0.5799930
\label{eq:plateau}
\end{equation}
On the right we also see a vertical 
jump, denoting the presence of a 
hyperbolic point. 

We are interested in the most robust tori (identified by their
frequencies) of the controlled system. They should be the last to
survive to the perturbation, and we want to target (one of) them with
the CAP of the next section. By the method of the Farey
tree \cite{kim_simultaneous_1986} we can take two rational numbers and
build a so called noble torus in between them, according to 
\cite{mackay_locally_1992} it is expected to be the most robust in a
local sense. The frequencies corresponding to the plateau are rational
numbers (because they correspond to resonant orbits): so we pose that
by the Farey tree, starting from the frequencies \eqref{eq:plateau},
we can find three robust tori of this dynamical system. Given two
rational numers \mbox{$n_1/d_1$} and
\mbox{$n_2/d_2$}, the ``most  noble'' 
number in between is 
\mbox{$(n_1+\sigma n_2)/(d_1+\sigma d_2)$},
being $\sigma = (\sqrt{5}+1)/2$, the 
so-called golden mean.

By trial and error, we have seen that 
a robust frequency is 
\begin{equation}
  \omega_D\,=\,\frac{43+\sigma 18}{74+\sigma 31}\,\approx\,0.580905,
  \label{eq:omegastar} 
\end{equation}
built by applying the Farey tree to
the couple $\omega_2,\,\omega_3$.

\section{CAP of the existence of a KAM torus}
\label{sec:cap}

In this section we describe the
scheme of a Computer Assisted Proof
(based on KAM theory) to prove the 
existence of an invariant 
torus for a Hamiltonian system.

In the constext of classical mechanics, 
KAM theory denotes a set of methods
to build invariant tori for perturbations 
of integrable Hamiltonian systems 
\cite{giorgilli_kolmogorov_1997}. 
To be more precise, it is possible to 
show that a dense set of invariant tori 
of an integrable Hamiltonian are only
deformed in presence of a (small) perturbation.
Let us remark that the control theory 
of M. Vittot (that we used in section 
\ref{sec:model}) is built in strong
analogy with KAM theory, as discussed 
in \cite{vittot_perturbation_2004}.

The KAM algorithm we are going to use was 
introduced in \cite{giorgilli_classical_1997} 
and \cite{giorgilli_kolmogorov_1997}; it 
was used for a CAP in 
\cite{celletti_improved_2000}, and then 
also in \cite{locatelli_invariant_2000}
and \cite{gabern_construction_2005}. 
It is based on classical series expansions for
the Hamiltonian function, with no use of the 
so-called quadratic method. By applying this 
algorithm to a Hamiltonian of the type 
\eqref{eq:hamiltonian2} for the model of the 
magnetic field, at the end of the $r$-th step 
of the algorithm, we get 
\begin{equation}
  \label{eq:hamiltonian_r}
  \hamcap^{(r)}\,
   =\,\omega\psi\,+\,\Pphi\,+\,\sum_{l=2}^{l_{\text{max}}}\sum_{s\geq 0} 
    \ham{l}^{(r,s)}
     \,+\,\sum_{s\geq r+1}\,\big(
      \delenda{0}^{(r,s)}
       \,+\,\delenda{1}^{(r,s)}
        \big)
\end{equation}
with $\ham{l}^{(r,s)}
\in\setpol{l}{sK}$, K can be defined 
as in section \ref{sec:model} as the
trigonometric degree of the perturbation 
$v$, or equivalently as the highest trigonometric
degree of the terms $\delenda{0,1}^{(0,s)}$.
By an infinite number of transformations, 
the initial Hamiltonian $\hamcap^{(0)}$ can be 
transformed to its Kolmogorov Normal Form
\begin{equation}
  \label{eq:KNF}
  \hamcap^{(\infty)}\,=\,\omega \psi\, +\,\Pphi\,+\,
   \sum_{l=2}^{l_\text{max}} \sum_{s\geq 0} 
    \ham{l}^{(\infty,s)}
\end{equation}
that has a manifestly invariant torus 
$\psi=0$ on which the motion takes 
place with frequency $\omega$.

In our CAP we begin by computing  
explicitly a finite number 
of KAM transformations, from $r=0$ 
to $r=R_I$. 
Then, we iterate estimates 
on upper bounds for the coefficients of 
the transformed Hamiltonians, from 
$r=R_I+1$ up to $r=\RII$.
In \cite{celletti_improved_2000} it
was discussed how a clever combination of 
the two techniques provides the best 
compromise between time-efficiency (the 
explicit computations are more demanding
in terms of computing time) and 
numerical precision (in the iteration of
the estimates we introduce a lot of 
approximations). Finally, we show that 
$\hamcap^{(\RII)}$
satisfies all the hypothesis of a KAM
theorem, so it can be conjugated by an 
analytical canonical transformation to 
the Kolmogov Normal Form
\eqref{eq:KNF} of $\hamcap^{(0)}$.

The present section just aims at giving a
well defined procedure leading to a complete
CAP, in such a way to perfectly fit with
the present context. Therefore, in the
following the justifications for
both the formal algorithm and the scheme
of the estimates will not
be discussed at all. Their validity can be
checked by referring to the already mentioned
works dealing with CAPs that are designed in
such a way to apply the KAM theorem.

In all of our codes that we wrote to perform the CAP
we implemented validated numerics so as to ensure 
fully rigorous results. A detailed discussion 
of validated numerics was given in appendix~A
of~\cite{caracciolo_computer-assisted_2020}

The codes performing the CAP are freely available from the {\it
  Mendeley Data} repository. They are collected in a software
package\footnote{That software package can be freely downloaded from
  the web address {\tt http://dx.doi.org/10.17632/jdx22ysh2s.1}}, to
which hereafter we will refer as the {\it supplementary material} that
is related to this paper. In order to increase their readability, they
are rather well commented and the structure of the whole software
package should be easy to understand starting from the detailed
explanations included in the file \texttt{README.txt}. We implemented
the explicit KAM transformations in the program \firstProg and the
iteration of the estimates in the program \texttt{iteration+proof.c};
the latter also computes the parameters that appear in the statement
of the KAM theorem.  In the following, we will refer sometimes to the
names of the functions that are written in {\bf C} programming
language and are included in the files making part of the {\it
  supplementary material}.  These references will be done in order to
eventually defer the reader to the algorithmic counterpart of our
explanations, when this can be seen as convenient for a better
understanding.

\subsection{Description of the KAM algorithm}
\label{sec:algo}

At the beginning of the $r$-th step of
the algorithm the Hamiltonian looks like
\begin{equation}
  \label{eq:hamiltonian_r-1}
  \hamcap^{(r-1)}\,=\,\omega \psi\,+\,\Pphi\,+\,
    \sum_{l=2}^{l_{\text{max}}}\, \sum_{s\geq 0}\, 
    \ham{l}^{(r-1,s)}
    \,+\,\sum_{s\geq r}\,\Big(
    \delenda{0}^{(r-1,s)}\,+\, 
    \delenda{1}^{(r-1,s)}\Big) 
\end{equation}
where (again) 
\mbox{$\ham{l}^{(r-1,s)}
\in\setpol{l}{sK}$}.
In principle the Hamiltonian has an infinite 
number of coefficients, but in our code
they are truncated to a maximum trigonometric
degree\footnote{
  \label{note:truncation}
  However, at each step we compute the norms 
  of the terms whose trigonometric 
  degree is bigger than \maxtrigdeg. They 
  are also exported and will be used as an 
  input for the iteration of the estimates.
} of\footnote{
  The factor $2$ stems from the number of
  angular variables that are considered in
  the system.
} \maxtrigdeg. For the representation of the 
Taylor-Fourier coefficients we employed a 
technique based on the so-called ``indexing
functions'' which is described in detail
in \cite{giorgilli_methods_2011}. From equation
\eqref{eq:hamiltonian_r-1} (but also from
equation
\eqref{eq:hamiltonian_r}) it is evident that
the main part of the angular average applied
to the linear term of the 
Hamiltonian (i.e., $\omega \psi\,+\,\Pphi$) 
has a different role with respect
to the other summands. Indeed, KAM 
algorithms are usually designed to preserve the
frequency $\omega$. In our code, this value 
is stored apart from the rest of the Hamiltonian, 
so it cannot be affected by numerical
errors related to the performing of the Lie
series\footnote{Of course, it is still subject to 
  the interval arithmetics.
}.

To the Hamiltonian \eqref{eq:hamiltonian_r-1}
we apply three Lie series, according to the rule 
\begin{equation}
  e^{\lieDer{g}}\hamcap\,=\,\sum_{j=0}^\infty 
   \,\frac{1}{j!}\,{\lieDer{g}^j} \hamcap
\label{eq:lieSeries}
\end{equation}
that is implemented in our code by calling the
function $\mathtt{can\_transf\_ham}$. 
The number of operations performed by
this function is always finite because 
the functions are truncated so that the 
maximum trigonometric degree of the result 
does not exceed \maxtrigdeg.

The first generating function
\mbox{$X^{(r)}\in\setpol{0}{rK}$}
is given by
\begin{equation}
  \label{eq:X}
  X^{(r)}\,=-\Gamma_{\omega}\, 
    \delenda{0}^{(r-1,r)}.
\end{equation}
As in section \ref{sec:model},
the operator $\Gamma_{\omega}$
is well defined if the frequency
$\omega$ is 
\textit{non-resonant}\footnote{
  This hypothesis is satisfied if
  $\omega$ is a Diophantine number,
  as we often assume in KAM theory.
}.
Also, the operator $\Gamma_{\omega}$ 
is well defined on 
a zero-average function; the average of 
$\delenda{0}^{(r-1,r)}$ is a constant and
can be safely ignored.
The effect of the Lie transform 
generated by $\lieDer{X^{(r)}}$
is exactly to erase $\delenda{0}^{(r-1,r)}$.

The operator $\Gamma_{\omega}$ is implemented
in our code by the function 
\texttt{solve\_homol\_eq} 
which moves the Fourier coefficients of
$\delenda{0}^{(r-1,r)}$
from the Hamiltonian to the 
variable $\mathtt{target\_fn}$, and 
use them to compute the coefficients 
of $X^{(r)}$ (stored in the
variable $\mathtt{gen\_fn}$).

The second generating function is
$\xi^{(r)} \theta$, and the corresponding
transform is a translation of the action
$\psi$ by a quantity $\xi^{(r)}$. 
The effect of this second transform 
is to erase $\hamavg{1}{(r-1,r)}$,
where $\hamavg{l}{(r,s)}$ is the average 
of a generic function $\ham{l}^{(r,s)}$
with respect to the
angle variables $\theta$ and $\varphi$,
for all \mbox{$0\leq l\leq l_\text{max}$}
and non-negative indexes \mbox{$r,\,s$}.
To fix the ideas, by referring to the
decomposition described
in equation~\eqref{eq:gp}, $\langle g\rangle$
is nothing but the coefficient $g_{l,0,0}$
corresponding to the harmonic $(k_1=0,k_2=0)$
appearing in the Taylor-Fourier expansion of
$g\in\setpol{l}{s}\,$.
Let us define 
\begin{equation}
  C^{(r)}\,=\,\sum_{s=0}^r \partial_{\psi\psi}^2 \hamavg{2}{(r-1,s)}
  \label{eq:defC}
\end{equation}
and impose 
\begin{equation}
  C^{(r)}\,\xi^{(r)}\,+\,\hamavg{1}{(r-1,r)}\,=\,0.
  \label{eq:xi}
\end{equation}

Equation \eqref{eq:xi}
is solved in our code by the function
$\mathtt{solve\_transl\_eq}$. For a 
$1\tfrac{1}{2}$d.o.f. system, $C$ is a scalar,
but in general it is a matrix
which needs to be inverted; this operation 
is performed by the function 
$\mathtt{calc\_inv\_matr\_C}$,
after removing the average of the
quadratic term of the Hamiltonian
through the function 
\texttt{move\_average\_to}.
Finally we need a function
\texttt{readd\_lie\_series\_csiq\_average\_terms}
to add back the average of the quadratic 
terms and the constant term that should
have been produced by the
action of the second Lie series on 
the average of the quadratric term, 
that instead we removed.

Sometimes the first two generating
functions are grouped into\footnote{
  This is expecially meaningful when the two trasforms
  are defined so as to commute; unfortunately, this is
  not our case, because to compute $\xi^{(r)}$ 
  we also need the coefficient
  $\ham{2}^{(r-1,r)}$ which is affected by the action of
  $X^{(r)}$. Nevertheless we will use this type
  of notation because it is useful for the
  next section. 
}
\mbox{$\chi_1^{(r)}\,=\,X^{(r)}
+\xi^{(r)}\theta$}.
Then, we introduce the 
intermediate Hamiltonian\footnote{
  \label{note:finite}
  Note that, by the action of 
  $e^{\lieDer{\chi_1^{(r)}}}$ on a Hamiltonian
  with a finite expansion in all variables, we get 
  another Hamiltonian with finite expansion. 
  In fact, as 
  \begin{equation}
  \lieDer{\chi_1^{(r)}} g\,=\,\big(\partial_\theta X^{(r)}(\theta,\varphi)\,+\,\xi^{(r)}\big)\partial_\psi g
  \notag 
  \end{equation}
  when we act with $e^{\lieDer{\chi_1^{(r)}}}$ on a 
  Hamiltonian like \eqref{eq:hamiltonian2},
  the Lie series can be applied $l_\text{max}$
  times at most. In particular, if 
  the trigonometric degree of the 
  original Hamiltonian is $R_h K$ and that
  of the generating function is $R_g K$ , the 
  transformed Hamiltonian has trigonometric
  degree $(R_h+l_\text{max}R_g)K$. 
}
\begin{equation}
 \begin{aligned}
  \hat{\hamcap}^{(r)}\,&=\,
    e^{\lieDer{\chi_1^{(r)}}}\hamcap^{(r-1)}\\
    &=\,\omega \psi\,+\,\Pphi\,+
    \sum_{l=2}^{l_{\text{max}}}\,\sum_{s\geq 0}
    \hatham{l}^{(r,s)}
    \,+\,\hatdel{1}^{(r,r)}
    \,+\,\sum_{s\geq r+1}\,\Big(
    \hatdel{0}^{(r,s)}
    +\,\hatdel{1}^{(r,s)}
    \Big)\,, \\
    \hatdel{1}^{(r,r)}\,&=\,
    \delenda{1}^{(r-1,r)}
     \,-\,\hamavg{1}{(r-1,r)}\, + \,
      \big(\partial_\psi \hamavg{2}{(0)}\big)
       \big( \partial_\theta X^{(r)} \big) .
 \end{aligned}
 \label{eq:hatH}
\end{equation}
After an inspection of equation 
\eqref{eq:hatH} we realize that we want 
to erase the term $\hatdel{1}^{(r,r)}$.
The latter is zero-average thanks to the
translation by $\xi^{(r)}$, so we 
can define the third generating 
function \mbox{$\chi_2^{(r)}\in\setpol{1}{rK}$} 
by another homological equation,
\begin{equation}
 \label{eq:chi2} 
 \chi_2^{(r)}\,=\,
  -\Gamma_{\omega}\,\hatdel{1}^{(r,r)}.
\end{equation}
The Lie series generated by $\chi_2^{(r)}$
maps $\hat{\hamcap}^{(r)}$ to the
the Hamiltonian \mbox{$\hamcap^{(r)}$}
of equation \eqref{eq:hamiltonian_r}.
To solve this second homological
equation we invoke again the function 
\texttt{solve\_homol\_eq}. 
In this second case, we still need 
the term $\hatdel{1}^{(r,r)}$
after its deletion. Indeed, we have 
\begin{multline}
  e^{\lieDer{\chi_2^{(r)}}}\Big(\omega\psi\,
   +\,\Pphi\,+\,\hatdel{1}^{(r,r)}\Big)\,=\\
  =\,\omega\psi\,+\,\Pphi\,+
  \lieDer{\chi_2^{(r)}}(\omega\psi+\Pphi)\,+\,\hatdel{1}^{(r,r)} 
  \,+\,\sum_{j\geq 1}\frac{1}{(j+1)!}
  \lieDer{\chi_2^{(r)}}^j\Big(\lieDer{\chi_2^{(r)}}(\omega\psi+\Pphi)\,
    +\,(j+1)\hatdel{1}^{(r,r)} \Big)\,=\\=\,
  \omega\psi\,+\,\Pphi\,+\,\sum_{j\geq 1}\frac{j}{(j+1)!}
  \lieDer{\chi_2^{(r)}}^j\hatdel{1}^{(r,r)}.
\label{eq:why_target}
\end{multline}
This means that we need $\hatdel{1}^{(r,r)}$
to transform the expansion of terms
appearing in $\hat{\hamcap}^{(r)}$ to that
corresponding to $\hamcap^{(r)}$. This is why
the function \texttt{solve\_homol\_eq}
copied $\hatdel{1}^{(r,r)}$ to 
the variable $\mathtt{target\_fn}$. 
The coefficients of the last term
of equation \eqref{eq:why_target} 
are computed by the function 
$\mathtt{readd\_lie\_series\_chi2\_target\_term}$.

The terms appearing into the expansions 
of $\hat{\hamcap}^{(r)}$ and 
$\hamcap^{(r)}$ 
can be computed according to the 
formulae (28)-(31) of\footnote{
  To translate to our notation 
  it is sufficient to make the
  replacements
  \begin{eqnarray}
    f_l^{(r,s)}\mapsto \ham{l}^{(r,s)}\,,\quad l\geq 0\,,\, s\geq r+1 \\
    h_l^{(s)}\mapsto \ham{l}^{(r,s)}\,,\quad l\geq 2\,,\, 0\leq s\leq r
  \label{eq:conversion}
  \end{eqnarray}
  and two analogous rules for 
  hatted quantities.
} \cite{locatelli_invariant_2000}.

We recall here that all the computations
are performed with interval arithmetics.
The machine represents 
real numbers by floating point numbers,
which however do not have the same properties
of reals. As 
a consequence, to make numerical results 
rigorous, it is possible to replace a
real number by \emph{two} floating point
numbers, representing the upper and lower
bound of the original real number in its 
floating point representation. For this 
purpose, in our code we employed 
an aptly defined type called 
\texttt{INTERVAL}. Evidently, all the 
elementary operations have to be redefined
for intervals; the principle to do so is
explained in appendix A of 
\cite{celletti_improved_2000} (short version) 
and also in appendix A of 
\cite{caracciolo_computer-assisted_2020}
(long version). In our code, the 
implementation of the operations for 
intervals is found in the file
\texttt{int\_arit\_lib.c}.

\subsection{Iteration of the Estimates}
\label{sec:iteration}

Given $g\in\setpol{l}{s}$ (see 
equation \eqref{eq:gp}) we introduce the norm
\begin{equation}
  \| g\|\,=\,\sum_{\max(|k_1|,|k_2|) \leq s} |g_{l,k_1,k_2}|.
  \label{eq:defNorm}
\end{equation}

So we introduce an upper bound for the 
norm of any function that we computed
explicitly in the preceding section.
First we define,
for any \mbox{$1\leq r\leq R_I$},  
the four constants\footnote{
  If the system was not degenerate in 
  the variable $\Pphi$, we 
  would replace 
  \mbox{$\partial_\theta\mapsto\max\big(\partial_{\theta},\partial_{\varphi}\big)$},
  \mbox{$|\xi^{(r)}|\mapsto \max\big(|\xi_1^{(r)}|,|\xi_2^{(r)}|\big)$},
  and finally 
  \mbox{$\partial_\psi\mapsto\max\big(\partial_{\psi}\,;\,\partial_{\Pphi}\big)$}.
}
\begin{equation}
  \genbound{1}{1}{r}\,=\,\big\| \partial_\theta X^{(r)} \big\|
  \,,\quad
  \genbound{1}{2}{r}\,=\,|\xi^{(r)}|
  \,,\quad
  \genbound{2}{1}{r}\,=\,\big\| \partial_\theta \chi_2^{(r)} \big\|
  \,,\quad
  \genbound{2}{2}{r}\,=\,\big\| \partial_\psi \chi_2^{(r)} \big\|.
  \label{eq:gs}
\end{equation}

Second, we need sets of constants to 
bound from above the terms appearing
in the expansions of 
the Hamiltonians $\hamcap^{(r)}$
and $\hat{\hamcap}^{(r)}$; for all 
$1\leq r\leq \RII$ we call them
\mbox{$\hbound{l}{r}{s}$},
\mbox{$\hathbound{l}{r}{s}$},
and \mbox{$\avgbound{r}{s}$},
and we define them so that, for all 
$1 \leq r \leq R_I$,
\begin{equation}
  \begin{aligned}
    &\vbound{l}{r}{s}\,=\,\|\delenda{l}^{(r,s)}\|
    \,,\quad
    \hatvbound{l}{r}{s}\,=\,\|\hatdel{l}^{(r,s)}\| 
    \,,\qquad l=0,1,\,r\leq s\leq R_I \,,\\
    & \hbound{l}{r}{s}\,=\,\|\ham{l}^{(r,s)} \|
    \,,\quad
    \hathbound{l}{r}{s}\,=\,\|\hatham{l}^{(r,s)}\| 
    \,,\qquad 2\leq l\leq l_{\text{max}},\, 0\leq s\leq R_I.
  \end{aligned}
  \label{eq:hv}
\end{equation}

The values of the right hand side of 
equation \eqref{eq:gs} and \eqref{eq:hv}
are computed by the program
\firstProg and written on
an output file  
\texttt{bounds\_expl\_expans.bin},
for all $1\leq r\leq R_I$. 
The values written on this file 
are read in the program 
\texttt{iteration+proof.c} by
the functions
\texttt{read\_gen\_fun\_bounds} and
\texttt{read\_ham\_bounds}.

We have now to face the problem of providing upper bounds for all
those terms that appear in the expansions of the Hamiltonians defined
by the normalization algorithm and are of type $\setpol{l}{sK}$ with
$0\leq l\leq l_{\text{max}}$ and $R_{I}< s \leq \RII\,$. Their
representations in Taylor-Fourier series is still finite, but they
include so many terms that such functions are hard to be explicitly
calculated. Nevertheless, an efficient scheme of iterative estimates
can be settled also for this kind of functions so to greatly improve
the final accuracy of the results provided by the CAP, when $\RII\gg
R_{I}\gg 1$, as it has been widely discussed
in~\cite{celletti_improved_2000}. For such a purpose, we have to
properly extend the definition of the multi-indicial finite sequencies
$\vbound{\cdot}{\cdot}{\cdot}$ and $\hatvbound{\cdot}{\cdot}{\cdot}$
so that the following inequalities hold true for all 
$1 \leq r \leq \RII\,$:
\begin{equation}
  \begin{aligned}
    &\|\delenda{l}^{(r,s)}\|\,\leq\,\vbound{l}{r}{s}
    \,,\quad
    \|\hatdel{l}^{(r,s)}\|\,\leq\,\hatvbound{l}{r}{s} 
    \,,\qquad l=0,1\,,\,r\leq s\,,\, R_I < s \leq \RII \,,\\
    & \|\ham{l}^{(r,s)} \|\,\leq\,\hbound{l}{r}{s}
    \,,\quad
    \|\hatham{l}^{(r,s)}\|\,\leq\,\hathbound{l}{r}{s} 
    \,,\qquad 2\leq l\leq l_{\text{max}}\,,\, 0\leq s\leq \RII.
  \end{aligned}
  \label{eq:maggiorazioniulteriori}
\end{equation}
Moreover, it is also convenient to define another multi-indicial
set of positive values $\avgbound{\cdot}{\cdot}$ so that
\begin{equation}
  \|< \hatdel{1}^{(r,s)} >\|\,\le\,\avgbound{r}{s}
  \,,\quad R_I\leq r\leq s\leq \RII .
\end{equation}
These three finite sequencies are fully defined according to the
prescriptions given in the following two
sections~\ref{sss:boundgenfun} and~\ref{sss:stimeiterative}.

\subsubsection{Bounds on the generating functions}
\label{sss:boundgenfun}

To grant the solvability of the homological
equation at each step, we assume that there
exists, for each value of the step index $r$, 
a positive constant $\alpha_r$ such that
\begin{equation}
  \big| k_1\omega + k_2\big|\,\geq\,\alpha_r
  \label{eq:dioph_2}
\end{equation}
and which can be simply computed by setting 
\begin{equation}
  \alpha_r\,=\,\min_{0<\max(|k_1|,|k_2|)\leq rK}
    \Big(\big|k_1 \omega + k_2\big|\Big).
  \label{eq:alpha_r}
\end{equation}
In both \eqref{eq:dioph_2} and \eqref{eq:alpha_r}
$k_1$ and $k_2$ are integers numbers that
cannot be simoultaneously zero.
We recall that the constant $K$ 
was defined in section 
\ref{sec:model} as the trigonometric 
degree of the perturbation $v$.
Note that the existence of a positive
constant $\alpha_r$ for $1\leq r\leq R_I$ is 
granted by the fact that the explicit
computation can be performed. We need to 
compute its value only for 
$R_I < r\leq \RII$.

From equation \eqref{eq:X}, we find
\begin{equation}
  \|X^{(r)}\|\,\leq\,\frac{1}{\alpha_r}\,\|\delenda{0}^{(r-1,r)}\|.
\label{eq:boundX}
\end{equation}
Therefore, by analogy with~\eqref{eq:gs} we define, for
all $R_I< r\leq \RII$,
\begin{equation}
  \genbound{1}{1}{r}\,=\,\frac{ rK }{ \alpha_r } \vbound{0}{r-1}{r} .
\label{eq:defg11}
\end{equation}
In a similar way, from equation
\eqref{eq:chi2} we can deduce 
\begin{equation}
  \|\chi_2^{(r)}\|\,\leq\,\frac{1}{ \alpha_r }\| \delenda{1}^{(r-1,r)} \| 
   \,\implies\,
    \genbound{2}{1}{r}\, =\,
     \frac{ rK \hatvbound{1}{r-1}{r} }{ \alpha_r }
      \,,\quad
       \genbound{2}{2}{r}\, =\,
        \frac{ \hatvbound{1}{r-1}{r} }{ \alpha_r }.
\label{eq:boundchi2}
\end{equation}
These equalities hold, again, for
$R_I< r\leq \RII$.

Finally, let us assume that there 
exists a 
positive constant $m^{(r)}$ such that 
$m^{(r)}\leq |C^{(r)}|$, then
\begin{equation}
  |\xi^{(r)}|\,\leq\,(m^{(r)})^{-1}\,\| \hamavg{1}{(r-1,r)} \|
    \implies \genbound{1}{2}{r}\,=\,(m^{(r)})^{-1}\,\avgbound{r}{r} .
\label{eq:boundxi}
\end{equation}
The direct computation of $m^{(R_I)}$ is 
obvious in our $1\tfrac{1}{2}$d.o.f.
example, and it could be also computed
in dimension $n>1$ as the eigenvalue
of the Hessian of $\hamcap^{(R_I)}$ 
with the lowest modulus.
Moreover, it is easy to show that
this constant can be updated
recursively according to the following
rule
\begin{equation}
  m^{(r)}\,=\,m^{(r-1)}\,-\,2\,\hbound{2}{r-1}{r}\,,
    \quad \forall\, R_I < r\leq \RII.
  \label{eq:mrec}
\end{equation}

\subsubsection{Bounds on the transformed Hamiltonian}
\label{sss:stimeiterative}

Given the bounds on the generating 
functions, we can deduce those on 
the terms appearing in the expansion of
a transformed Hamiltonian $H^{(r)}$ starting
from that corresponding to $H^{(r-1)}$ for
all $1\le r\le\RII\,$. In the case of the
transformation generated by $\chi_1^{(r)}$,
\begin{equation}
 \Big\| \frac{1}{j!}\,\lieDer{\chi_1^{(r)}}^j g \Big\| \,\leq\,\binom{l}{j}\Big(
  \big\|\partial_\theta X^{(r)} \big\| \, +\, \big| \xi^{(r)} \big| \Big)^j\, \| g \| .
 \label{eq:boundsLie1}
\end{equation}
If we assume that \mbox{$g\in\setpol{l}{sK}$},
then
\mbox{$\lieDer{\chi_1^{(r)}}^j g\in\setpol{l-j}{(s+rj)K}$}.
So we can update the upper bounds on
the terms appearing in the expansion of
the Hamiltonian $H^{(r-1)}$ (that appear in
formula~\eqref{eq:hamiltonian_r-1} and
are assumed to be initialized when
$r-1=0$) according to the
following rules\footnote{
  We introduce the symbol $\pluseq$
  (familiar to the C programmers)
  so that $a\pluseq b$
  is equivalent to a {\it redefinition}
  of the value of $a$ in such a way
  that $a\mapsto a+b$.
}: $\hatvbound{l}{r}{s}=\hbound{l}{r-1}{s}$ $\forall\, l,\,s$,
$\hatvbound{0}{r}{r}=0$ and
\begin{equation}
 \begin{gathered}
  \hathbound{l-j}{r}{s+jr}\,\pluseq\,\binom{l}{j}
   \Big(\genbound{1}{1}{r}\,+\,\genbound{1}{2}{r}\Big)^j
   \hbound{l}{s}{s}\,,\qquad \forall\, l\geq 2,\,0\leq s< r,\,1\leq j\leq l,
   \,R_I<s+jr\leq\RII\,,\\
  \hatvbound{l-j}{r}{s+jr}\,\pluseq\,\binom{l}{j}
   \Big(\genbound{1}{1}{r}\,+\,\genbound{1}{2}{r}\Big)^j
    \hbound{l}{r-1}{s}\,,\qquad \forall\, l\geq 1,\,s\ge r,\,1\leq j\leq l,
    \,R_I<s+jr\leq\RII. \\
 \end{gathered}
\label{eq:hatcoef}
\end{equation}

In the case of the trasformation generated by 
$\chi_2^{(r)}$, the following inequality holds
true for any $g\in\setpol{l}{sK}\,$:
\begin{equation}
 \Big\| \frac{1}{j!}\,\lieDer{\chi_2^{(r)}}^j g \Big\| \,\leq\,
  \frac{1}{j!}\,\prod_{i=0}^{j-1} \Big[ l\big\|\partial_\theta \chi_2^{(r)}\big\|
   \,+\, \big( (s+ir)K \big) \big\|\partial_\psi \chi_2^{(r)} \big\|  \Big] \,
    \| g \|.
 \label{eq:boundsLie2}
\end{equation}
Since 
\mbox{$\lieDer{\chi_2^{(r)}}^j g\in\setpol{l}{(s+rj)K}$},
a second set of rules easily follows:
$\hbound{l}{r}{s}=\hatvbound{l}{r}{s}$ $\forall\, l,\,s$,
$\hbound{1}{r}{r}=0$ and
\begin{equation}
 \begin{gathered}
  \hbound{l}{r}{s+jr}\,\pluseq\,
  \frac{1}{j!}\,\prod_{i=0}^{j-1} \Big[ l\genbound{2}{1}{r} \,+\,
    \big( (s+ir)K \big) \genbound{2}{2}{r} \Big] \hathbound{l}{s}{s}\,,
     \qquad \forall\, l\geq 2,\, 0\leq s\leq r,\, j\geq 1,\,R_I<s+jr\leq\RII\,, \\
   \vbound{l}{r}{s+jr}\,\pluseq\,
   \frac{1}{j!}\,\prod_{i=1}^j \Big[ l\genbound{2}{1}{r}\,+\,
     \big( (s+ir)K \big) \genbound{2}{2}{r} \Big] \hathbound{l}{s}{s}\,,
     \qquad \forall\, l\geq 0,\, s>r,\, j\geq 1,\,\,R_I<s+jr\leq\RII\,,
   \\ 
   \vbound{1}{r}{(j+1)r}\,\pluseq\,
   \frac{1}{j!}\,\prod_{i=1}^j \Big[ \genbound{2}{1}{r} \,+\,
     \big( (i+1)rK \big)
    \genbound{2}{2}{r} \Big] \hathbound{1}{r}{r}
     \,,\qquad \forall\, j\geq 1,\,R_I<(j+1)r\leq\RII.
\end{gathered}
\label{eq:newcoef}
\end{equation}

For what concerns the angular average of the perturbing terms that are
linear in $\psi$, we need to refine the iterative definitions with
respect to those provided for $\vbound{1}{r}{s}$.  This can be made by
using sharper evaluations preventing a deterioration effect that
would be due to the fact that the majorant $\genbound{1}{2}{r}$ would
contribute to the definition of $\genbound{1}{2}{s}$ for any
$s=r+1,\,\ldots\,2r-1$. Indeed, a so strict recursion is artificious
and can be avoided by slightly adapting the scheme of estimates in the
following suitable way. After having recalled the
equations~\eqref{eq:X}--\eqref{eq:why_target}, for every $r\ge R_I+1$
it is convenient to set
\begin{equation}
 \begin{gathered}
   \avgbound{r}{s}=\hbound{1}{r}{s}\,,
   \qquad {\rm if}\ r=R_I+1\ \forall\,r\leq s\leq\RII\,,
   \ {\rm else}\ \forall\,2r-2\leq s\leq\RII\,,
   \\
   \avgbound{r}{s+r}\,=\,\avgbound{r-1}{s+r}+2\genbound{1}{1}{r}
   \hbound{2}{r-1}{s}\,,\qquad \forall\,1\leq s< r,\,
   \,R_I<s+r\leq\RII. \\
 \end{gathered}
 \label{eq:def-avgbound}
\end{equation}
  
Formulae appearing in~\eqref{eq:hatcoef}
and~\eqref{eq:newcoef} are implemented 
in our program \texttt{iteration+proof.c} by the 
functions~$\mathtt{calc\_estimates\_ham\_chi1}$
and~$\mathtt{calc\_estimates\_ham\_chi2}$,
respectively.
Both are invoked for all 
$1\leq r \leq \RII$. Here, let us recall
that the constants $\hbound{l}{r}{s}$
need not to be computed by the recursive relations for 
$r,s\leq R_I$, because they are evaluated as norms
of functions, whose expansions are explicitly
computed according to the prescriptions
described in subsection~\ref{sec:algo}.
Formulae appearing in~\eqref{eq:def-avgbound}
are implemented by the 
function~$\mathtt{calc\_estimates\_ave\_ham\_chi1}$
that is invoked for all $R_I+1\leq r \leq \RII\,$.

\subsubsection{The recursive relations defining the parameters $\ecap_{r}$, $\zcap_{r}$ and $a_{r}$}

Of course, the expansions of the Hamiltonians contain an
infinite number of terms. Therefore, we need upper bounds
also for the norms of the functions $\ham{l}^{(r,s)}$ with
$s> \RII$. We stress that this kind of terms is neither
explicitly calculated according the algorithm described
in subsection~\ref{sec:algo} nor estimated as it has been
explained in the previous subsection~\ref{sss:stimeiterative}.
Therefore, we need to introduce three sequencies 
$(\ecap_r,\zcap_r,\acap_r)_{1\leq r\leq\RII}$, defined
so that\footnote{
  see for instance \cite{giorgilli_kolmogorov_1997},
  \cite{celletti_improved_2000},
  and \cite{locatelli_invariant_2000}.
}
at any step $r=1\,,\,\ldots\,,\,\RII$ of the KAM
algorithm the following inequalities are satisfied:
\begin{equation}
  \|\ham{l}^{(r,s)}\|\, \leq\, \ecap_r \,\zeta_r^l\, \acap_r^s\,,
        \quad \forall\, s>\RII\,,\,0\leq l\leq l_{\text{max}}. \\
  \label{eq:ineq}
\end{equation}

It is possible to derive recursive rules to 
compute these sequencies from the algorithm 
of section \ref{sec:algo} (see for instance
\cite{giorgilli_kolmogorov_1997},
\cite{celletti_improved_2000}). These rules,
which are not unique, 
involve also the norms of the generating functions. 
For $1\leq r\leq R_I$, these 
norms were computed explicitly. Therefore,
the bounds that we determine by the recursive
definition of
\mbox{$(\ecap_r,\zcap_r,\acap_r)_{r=R_I+1}^{\RII}$}. 
are likely to be
less sharp than those that we find 
when the explicit computations of the generating
functions are involved.

In the KAM algorithm of section \ref{sec:algo},
during the $r$-th normalization step we
transform $\hamcap^{(r-1)}$ into an
intermediate Hamiltonian 
$\hat{\hamcap}^{(r-1)}$ by the action of
the Lie series operator $e^{\lieDer{\chi_1^{(r)}}}$,
then  we transform $\hat{\hamcap}^{(r-1)}$ into 
$\hamcap^{(r)}$ by the action of 
$e^{\lieDer{\chi_2^{(r)}}}$. Accordingly, we choose
to update twice also the recursive
inequalities, at each step, in the following
way. Given $\ecap_{r-1}$, $\zcap_{r-1}$ and 
$\acap_{r-1}$, we compute
\begin{equation}
\begin{gathered}
  \hat{\ecap}_r\,=\,\ecap_{r-1} \bigg( 1 + 
   \frac{ \big(\, \genbound{1}{1}{r} + \genbound{1}{2}{r}\, \big) \zcap_{r-1} }
    { \acap_{r-1}^r } \bigg)^{l_{\text{max}}}
  \,,\quad
  \hat{\zcap}_r\,=\,\zcap_{r-1}
   \,,\quad
    \hat{\acap}_r\,=\,\acap_{r-1}\,,\qquad {\rm when}\ r\leq R_I\,, \\
  \hat{\ecap}_r\,=\,\ecap_{r-1}
  \,,\quad
  \hat{\zcap}_r\,=\,\zcap_{r-1}
   \,,\quad
   \hat{\acap}_r\,=\,\acap_{r-1}\bigg( 1 + l_{\text{max}}\,
   \frac{ \big(\, \genbound{1}{1}{r} + \genbound{1}{2}{r}\, \big) \zcap_{r-1} }
    { \acap_{r-1}^{r} } \bigg)^{\frac{1}{r}}\,,\qquad {\rm if}\ R_I< r\leq \RII\,.
\end{gathered}
\label{eq:rec_hat}
\end{equation}
and then
\begin{equation}
\begin{gathered}
  \ecap_r\,=\,\hat{\ecap}_r \,,\\
  \zcap_r\,=\,\hat{\zcap}_r
     \bigg(1+ \frac{\max\big(\, rK\genbound{2}{2}{r}\, ,\, \genbound{2}{1}{r}\,\big) }
          {{\hat\acap}_{r}^r} \bigg)\,,\quad
  \acap_r\,=\,\hat{\acap}_r\bigg(1+ \frac{
   \max\big(\, rK\genbound{2}{2}{r} \,,\, \genbound{2}{1}{r} \,\big)}
    {{\hat\acap}_{r}^r} \bigg)^{\tfrac{1}{r}}\,,\qquad {\rm when}\ r\leq R_I\,, \\
  \zcap_r\,=\,\hat{\zcap}_{r}\,,\quad        
  \acap_r\,=\,\hat{\acap}_r\bigg( 1+ \frac{(l_\text{max}+1)
   \max\big(\,rK\genbound{2}{2}{r}\,,\,\genbound{2}{1}{r}\,\big)}{\hat{\acap}_{r}^r} 
  \bigg)^{\tfrac{1}{r}}\,,
  \qquad {\rm if}\ R_I< r\leq \RII\,.
\end{gathered}
\label{eq:rec_new}
\end{equation}
Let us remark that we could keep $\zcap_r$ constant
at the price of introducing a factor $l_{\text{max}}$
in the rule for $\acap_r$, but this choice is
advantageous (because the contribution of
such a factor is made negligible) only when
the exponent $1/r$ is small.

The program \texttt{iteration+proof.c} also
invokes orderly the functions 
$\mathtt{update\_estimates\_param\_chi1}$
and $\mathtt{update\_estimates\_param\_chi2}$
to compute, respectively, the
values of $\hat{\ecap}_r,\,\hat{\zcap}_r,\,\hat{\acap}_r$
and of $\ecap_r,\,\zcap_r,\,\acap_r\,$,
for all $1 \leq r \leq \RII\,$.

To define the initial values $\acap_1$, 
$\zcap_1$ and $\ecap_1$ we exploit the 
fact that the original Hamiltonian
has a finite degree in both actions 
and angles; let us recall that they
are equal to ${l_\text{max}}$ and
$R_I K$, respectively.
The sequence $\{\acap_r\}_r$ is the most
threatening one to achieve convergence 
of the CAP; it is non-decreasing and we
want to initialize it to the smallest
value compatible with the inequalities
in formula~\eqref{eq:ineq}, when $r=1$.
Therefore, we set
\begin{equation}
  \acap_1\,=\,\max\big(\,K\genbound{2}{2}{1}\,,\,\genbound{2}{1}{1}\,\big)\,,
  \qquad \zcap_1\,=\,1\,,\qquad
  \ecap_1\,=\,\max_{\{0\leq l\leq l_\text{max}\,,\, 0\leq s\leq R_I+l_\text{max}\}} 
    \bigg( \frac{ \hathbound{l}{1}{s} }{\acap_1^s\zcap_1^l} \bigg).
\label{eq:init_rec}
\end{equation}
These values are computed by the 
function \texttt{initialize\_estimates\_param}.

\subsection{Statement of a KAM theorem}
\label{sec:thm}

The formulation of KAM theorem we are 
going to consider is the one reported
in~\cite{stefanelli_kolmogorovs_2012},
that is adapted to a slightly larger
class of dynamical models. Indeed, it
applies to systems including also a particular
type of dissipation that is proportional to
a constant factor $\eta$, but also the
conservative case with $\eta=0$ is covered.
Moreover, the framework of the proof of
that version of the KAM theorem is different
with respect to what is assumed here,
therefore, a little of effort for the
adaptation to the present context is mandatory.
For instance,
in~\cite{stefanelli_kolmogorovs_2012} the
\textit{weighted Fourier norms} are
considered. In order to properly introduce them,
let us first define the following complex sets:
$\mathbb{D}_{(\rho,\sigma)}\,=\,\mathbb{G}_\rho\times \mathbb{T}_{\sigma}^2\,$,
where
\begin{equation}
  \label{eq:DrhoSigma}
  \mathbb{G}_\rho\,=\,
    \big\{ z\in\mathbb{C}\,\colon\,
      |z|<\rho \big\}\ ,
  \qquad
  \mathbb{T}_{\sigma}^2\,=\,
   \big\{\,(\tilde{\theta},
    \tilde{\varphi})\in\mathbb{C}^2\,\colon\,
     \big(Re(\tilde{\theta}),Re(\tilde{\varphi})\big)\in\mathbb{T}^2,\,
       |Im(\tilde{\theta})|\leq\rho,\, 
        |Im(\tilde{\varphi})|\leq\rho\,\big\}\ ,
\end{equation}
being $\rho$ and $\sigma$ two positive real parameters.
On the set $\mathbb{D}_{(\rho,\sigma)}$ we
can define the functional norms  
\begin{equation}
 \|g\|_{(\rho,\sigma)}\,=\,
  \sum_{ (k_1,k_2) \in \mathbb{Z}_2/(0,0) }
   \sup_{\psi\in\mathbb{G}_\rho}|g_{l,k_1,k_2}|e^{(|k_1|+|k_2|)\sigma}.
\label{eq:defWFN}
\end{equation}
In order to concile our notations with those 
of \cite{stefanelli_kolmogorovs_2012} 
let us now introduce
\begin{equation}
 \delenda{0,1}^{(\RII)}\,=\,\sum_{s\geq \RII+1}
  \delenda{0,1}^{(\RII,s)}\,,\quad
 \ham{l}^{(\RII)}\,=\,\sum_{s\geq 0}
  \ham{l}^{(\RII,s)}\,,2\leq l\leq l_\text{max}\,.
\label{eq:defCoefThm}
\end{equation}

Theorem 3.1 of 
\cite{stefanelli_kolmogorovs_2012} 
can now be stated in the following way
to fit with the context of the present
approach.

\begin{theorem}
\label{thm:kam}
Consider a Hamiltonian
$H^{(\RII)}=\omega \psi +\Pphi+ \sum_{l\ge 0}\ham{l}^{(\RII)}$,
that is analytic on the domain
$\mathbb{D}_{(\rho,\sigma)}=
\mathbb{G}_\rho\times \mathbb{T}_{\sigma}^2$
defined in formula~\eqref{eq:DrhoSigma}.
Let us assume that 
\begin{itemize}
\item there exist positive numbers 
  $\rho,\sigma,\lambda,E$ such that
  \begin{equation}
    \big\| \delenda{0}^{(\RII)} \big\|_{(\rho,\sigma)}\,\leq
    \, \lambda E\,;\quad
    \big\| \delenda{1}^{(\RII)} \big\|_{(\rho,\sigma)}\,\leq
    \, \frac{\lambda E}{2}\,;\quad
    \big\| \ham{l}^{(\RII)} \big\|_{(\rho,\sigma)}\,\leq
    \, \frac{E}{2^l}\,,\ \forall\, l\ge 0 \ ;
  \label{eq:boundsThm}
  \end{equation}
\item there are two positive constants
  $\gamma$ and $\tau\geq 1$ such
  that\footnote{This is usually known
  as the \emph{Diophantine condition}.}
  $|k_1\omega+k_2|\ge \gamma/(|k_1|+|k_2|)^\tau$
  $\forall\, (k_1\,,\,k_2)\neq (0,0)$;
\item there exists a positive constant $m$
  such that $\big|\partial_{\psi\psi}^2\hamavg{2}{(\RII)}\big|\, \geq\, m$;
\item the parameter $\lambda$ is so small that
\begin{equation}
  \lambda < \lambda^*\,=\,\min\bigg( 
   \frac{1}{20^{25+9\tau}\,\bar{W}^2\,\bar{Z}^2}\,;
    \,9\Big(\frac{m\rho^2}{Ee^{\pi^2/3}}\Big)^2\bigg)
\label{eq:lambda*}
\end{equation}
where
\begin{eqnarray}
  \bar{A}\,=\,\frac{ E e^{\pi^2/3} }{\gamma}\Big(\frac{\tau}{e\sigma}\Big)^\tau \,,\qquad 
  \bar{B}\,=\,\Big( 1 + \frac{ \bar{A} }{ e\rho\sigma } \Big) \frac{ E e^{\pi^2/3} }{m\rho}\ ,  \\
  \bar{W}\,=\,\max\bigg( \frac{e^2}{\rho}\Big(\frac{\bar{A}}{e\sigma}\,+\,\bar{B}\Big)\,,\, 2 \bigg)  \,,\qquad 
  \bar{Z}\,=\,\max\bigg( \frac{2e\bar{A}}{\rho\sigma}\,,\,2 \bigg)\ .
\label{eq:paramThm}
\end{eqnarray}
\end{itemize}
Therefore, there exists a canonical analytical 
change of variables 
\mbox{$ {\cal C} \colon \mathbb{D}_{ \tfrac{1}{2}(\rho,\sigma) }
  \to \mathbb{D}_{ \tfrac{3}{4}(\rho,\sigma) }$}, 
transforming $H^{(\RII)}$ in the Kolmogorov 
Normal Form of type~\eqref{eq:KNF}.
\end{theorem}

In order to apply the statement above,
we need to choose four values of
the parameters $\sigma$, $\rho$,
$E$ and $\lambda$ so that the inequalities
in~\eqref{eq:boundsThm} are satisfied.
For this purpose, we have to reconsider
the upper bounds that we introduced
in the previous section~\ref{sec:iteration}.
Let us remark that 
\begin{equation}
 \big\|\ham{l}^{(\RII,s)}\big\|_{(\rho,\sigma)}\,\leq
  \big\|\ham{l}^{(\RII,s)}\big\|\,\rho^l\,(e^{2K\sigma})^s
\label{eq:trasl_bounds}
\end{equation}
(where the form of the last multiplying
coefficient is due to the fact that
when the norm of $\ham{l}^{(\RII,s)}$ is
evaluated according to~\eqref{eq:defWFN}
we have that $|k_1|+|k_2|\le 2sK$, being
the number of degrees of freedom equal
to $2$) and so 
\begin{eqnarray}
\label{eq:bound_0}
 &\big\|\delenda{0}^{(\RII)}\big\|_{(\rho,\sigma)}\,\leq\,
  \ecap_{\RII}\sum_{s\geq \RII+1} (a_{\RII} e^{2K\sigma})^s\,
   =\,\ecap_{\RII} 
   \frac{ (a_{\RII} e^{2K\sigma})^{(\RII+1)} }{1-a_{\RII} e^{2K\sigma}}\,, \\[0.2cm]
\label{eq:bound_1}
 &\big\|\delenda{1}^{(\RII)}\big\|_{(\rho,\sigma)}\,\leq\,
  \ecap_{\RII}\sum_{s\geq \RII+1}\, 
   (\rho \zcap_{\RII})\, (a_{\RII} e^{2K\sigma})^s\,
   =\,\ecap_{\RII} \, (\rho \zcap_{\RII})\,
   \frac{ (a_{\RII} e^{2K\sigma})^{(\RII+1)} }{1-a_{\RII} e^{2K\sigma}}\,, \\[0.2cm]
\label{eq:bound_l}
  &\big\|\ham{l}^{(\RII)}\big\|_{(\rho,\sigma)}\,\leq\,
  \sum_{s=0}^{\RII}\,\hbound{l}{\RII}{s}\rho^l (e^{2K\sigma})^s\,+\,
   \sum_{s\geq \RII+1}\, \ecap_{\RII}\,(\rho \zcap_{\RII})^l\, 
    (a_{\RII} e^{2K\sigma})^s
          \,,\quad 2\leq l\leq l_\text{max}\,.
\end{eqnarray}
By comparison of the
formulae~\eqref{eq:bound_0}--\eqref{eq:bound_l}
and~\eqref{eq:boundsThm}, it looks quite
natural to  set $\rho\,=\,(2\zcap_{\RII})^{-1}$.
Moreover, the summation of the geometric
series 
\begin{equation}
  \sum_{s\geq \RII+1} (a_{\RII} e^{2K\sigma} )^s\,=\,
   \frac{ (a_{\RII} e^{2K\sigma})^{\RII+1} }{1-a_{\RII} e^{2K\sigma}} 
\label{eq:sumgeom}
\end{equation}
is well defined provided that
\begin{equation}
  \sigma\,<\, -\frac{1}{2K}\,\log( a_{\RII} )\,.
\label{eq:condConvGeom}
\end{equation} 
In order to define more conveniently the
relation among $\sigma$ 
and the other parameters, we try to minimize 
the value of $\lambda/\lambda^*$ with respect
to $\sigma$; with minor approximations we end 
up with \mbox{$\sigma=4/(K\RII)$}.
Finally, on the r.h.s. of equation
\eqref{eq:bound_l} we can group a factor
\mbox{$(\rho\zcap_{\RII})^l=2^{-l}$}, and we 
are naturally lead to define
\begin{equation}
 E\,=\,\max_{2\leq l\leq l_\text{max}}\left(\zcap_{\RII}^{-l}\, 
  \sum_{s=0}^{\RII} 
   \hbound{l}{\RII}{s} (e^{2K\sigma})^s
    \,+\,\ecap_{\RII}\,\frac{ (a_{\RII} e^{2K\sigma})^{(\RII+1)} }
     { 1-a_{\RII} e^{2K\sigma} }\right)
 \label{eq:E}
\end{equation}
and
\begin{equation}
 \lambda\,=\,\frac{\ecap_{\RII}}{E}\,
  \frac{ (a_{\RII} e^{2K\sigma})^{\RII+1} }
  {1-a_{\RII} e^{2K\sigma}}\,. 
 \label{eq:lambda}
\end{equation}

The parameter $\gamma$ (that is
entering the Diophantine
inequality reported in the
statement of theorem~\ref{thm:kam}) is
also computed in the program 
\texttt{iteration+proof.c}, by the 
function \texttt{calc\_gamma}.
The parameter $m$ that is also appearing
in the hypothesis of theorem 
\ref{thm:kam}, can be identified
with $m^{(\RII)}$ computed 
according to equation \eqref{eq:mrec}. 
The value of $m^{(R_I)}$ is computed
by the program \texttt{expl\_transf.c}
with the name \texttt{coef\_nondeg\_matr\_C},
and written to the output file
\texttt{bounds\_expl\_expans}; then
it is read by the program
\texttt{iteration+proof.c}
and stored in the variable
\texttt{m\_nondeg} and updated for all
\mbox{$R_I+1\leq r \leq \RII$}.
Moreover, the values of the four
parameters $\sigma$, $\rho$,
$E$ and $\lambda$
are all computed by the function
\texttt{calc\_param\_proof} in
the program \texttt{iteration+proof.c}.

Therefore, at the end of its execution,
the program \texttt{iteration+proof.c} can
check if theorem~\ref{thm:kam} applies
to the Hamiltonian $H^{(\RII)}$. This is
done by simply comparing the
small parameter $\lambda$ with its
threshold value $\lambda^*$ that is
computed by using the definitions
appearing in
formulae~\eqref{eq:lambda*}--\eqref{eq:paramThm}.

\section{Application to the problem of the controlled magnetic field}
\label{sec:appo}

Here we want to apply the CAP (described
in the previous section \ref{sec:cap}) to
the problem of the controlled magnetic
field (described in section \ref{sec:model}). 
The controlled Hamiltonian is 
$\tilde{H}+f$, see equations
\eqref{eq:hamiltonian2}, and
\eqref{eq:control2}. 

In the qualitative analysis of section 
\ref{sec:afma} we have shown that 
we expect there is an invariant torus
having frequency $\omega=\omega_D$
(see equation \eqref{eq:omegastar}),
which is not the frequency corresponding to
the torus that is labelled by $\psi=0$
and is invariant with respect to the
flow of an integrable approximation of
$\tilde{H}+f$.
It is convenient that such a property
will be satisfied by $\hamcap^{(0)}$, i.e.,
the initial Hamiltonian that is
provided as an input for the whole CAP,
which is based on the KAM algorithm
(see section~\ref{sec:algo}).
The second type of
Lie series, that we apply at each 
step of the KAM algorithm, performs
a translation of the action and so
it yields also
a change of the frequency for
the torus corresponding to $\psi=0$
in the approximation provided by the
Kolmogorov normal form. Therefore, we 
decided to define $\hamcap^{(0)}$ 
by applying to $\tilde{H}+f$
a normalization step of the KAM
algorithm as follows
\begin{equation}
 H_c\,=\,e^{\lieDer{\chi_2}} 
  e^{\lieDer{\xi_\text{init}\theta}}
   e^{\lieDer{X}}\big(\tilde{H}+f\big)\,.
 \label{eq:Hc}
\end{equation}
The value of $\xi_\text{init}$ is defined 
to bring the frequency we can associate
to the torus corresponding to $\psi=0$
for the integrable approximation of
\mbox{$H_c$}
as close as possible to $\omega_D$, as we 
explain in the next
subsection~\ref{sec:xistar}.
Moreover, the generating function
$X$ is determined according to
formula~\eqref{eq:defGamma}, while
also $\chi_2$ is given by an equation
that is similar to that reported
in~\eqref{eq:defGamma}, but it is
aiming to remove the terms that
are linear in $\psi$ and with
Fourier harmonics
\mbox{$(k_1\,,\,k_2)$} such that
\mbox{$0<\max(|k_1|,|k_2|)\leq K$}.

\begin{figure}
\centering
\includegraphics[width=0.7\textwidth]{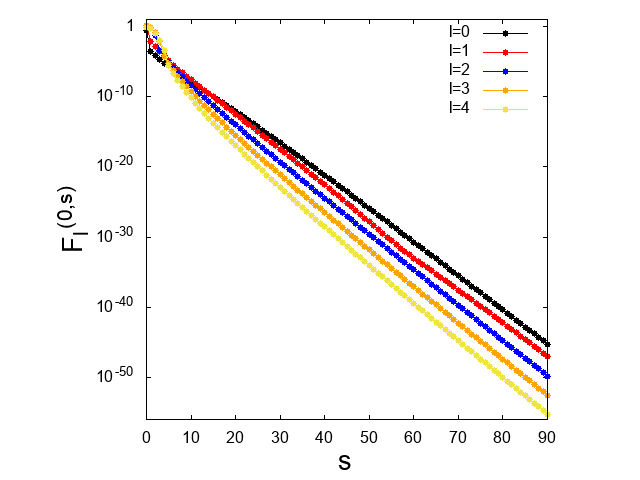}
\caption{
  \label{fig:norms_vs_s}
  Decay of the norms of the terms
  appearing in the
  expansion~\eqref{eq:Hc_as_H0}
  of $H^{(0)}$,
  that is defined in
  formulae~\eqref{eq:def-H0_parte1}--\eqref{eq:def-H0_parte2}. 
  }
\end{figure}

We can now introduce
\begin{equation}
  H^{(0)}\,=\,\omega_D\psi\,+\,\Pphi\,+\,
   \sum_{l=0}^{l_\text{max}} 
   \sum_{s=0}^{R_I} F_l^{(0,s)}
\label{eq:Hc_as_H0}
\end{equation}
where $F_l^{(0,s)}\in\setpol{l}{sK}$.
$H^{(0)}$ is meant to be a truncation
of the controlled Hamiltonian $H_c$
in the following sense. Let us write
the Taylor-Fourier expansion for both
of them
\begin{equation}
\label{eq:def-H0_parte1}
H_c=\sum_{l=0}^{l_\text{max}}\,\sum_{ (k_1\,,\,k_2)\in \mathbb{Z}^2}
      \,c_{l,k_1,k_2}\,\psi^l\,e^{ik_1\theta+ik_2\varphi} \ ,
\qquad
H^{(0)}=\sum_{l=0}^{l_\text{max}}\,\sum_{ (k_1\,,\,k_2)\in \mathbb{Z}^2}
      \,d_{l,k_1,k_2}\,\psi^l\,e^{ik_1\theta+ik_2\varphi} \ .
\end{equation}
Then
\begin{equation}
\label{eq:def-H0_parte2}
d_{l,k_1,k_2} =
\left\{
\vcenter{\openup1\jot 
\halign{
$\hfil\displaystyle{#}\hfil$\quad&$\displaystyle{#}\hfil$\cr
c_{l,k_1,k_2} &{\rm if}\ \max(|k_1|,|k_2|)\leq R_IK
\cr
0 &{\rm otherwise}
\cr
}}
\right .
\end{equation}
so that we have truncated the
Fourier coefficients in such a way
that the maximum trigonometric degree
of the terms appearing in $H^{(0)}$ is
not larger than
\maxtrigdeg\ (see section \ref{sec:algo})
before its expansion is passed to the 
program \texttt{expl\_transf.c}. 
In figure \ref{fig:norms_vs_s} we plot the
magnitude of the norms of the terms
$\ham{l}^{(0,s)}$ as a function of $s$: we see
that the truncation of the Fourier harmonics
at \mbox{\maxtrigdeg} is not a big issue
when $R_I$ is large enough.
In fact, for $s>90$ the magnitude of 
$\|\ham{l}^{(0,s)}\|$ is even smaller than
$10^{-45}$ $\forall\,l=0\,,\,\ldots\,,\,l_\text{max}\,$.

The initial Hamiltonian $H^{(0)}$
is preliminarly calculated by using
{\it X$\rho$\'o$\nu o\varsigma$}, that is
a software package especially designed
for doing computer algebra manipulations
into the framework of Hamiltonian
perturbation theory
(see~\cite{giorgilli_methods_2011} for
an introduction to its main concepts).
Such a code allows to represent the
complete (and finite) expansion in
Taylor-Fourier series of $H^{(0)}$,
in such a way to store it in the input
(ASCII) file \texttt{initial\_Ham0.inp},
that is read by the function
\texttt{def\_init\_ham}, which makes part
of the program \texttt{expl\_transf.c}.
A couple of intervals including the
components of the angular velocity vector
$\big(\omega_D\,,\,1\big)$ are written
in another input file, namely
\texttt{freq\_intervals.inp}.
Also for the sake of the usability of the
codes making part of the {\it supplementary
material}, we think it is preferable
to not include there any software package
that is depending on
{\it X$\rho$\'o$\nu o\varsigma$}.
Therefore, in that {\it material} we have
been forced to include the input files,
that are \texttt{initial\_Ham0.inp}
and \texttt{freq\_intervals.inp}.
More precisely, they are stored in a
subfolder, while in the main folder
there are two different versions of
those input files, that refer to an
example of a slightly perturbed
forced pendulum model, which is studied
in~\cite{celletti_improved_2000}.
This has been made for the sake of
the (re)usability of the package
provided in the {\it supplementary
material}. Since the Taylor-Fourier
representation of that Hamiltonian
system is extremely shorter with respect
to $H^{(0)}$, which is defined as discussed
above, our choice should allow to easily
understand how the input files must be
structured and modified in order to
design further different applications.

Our main result is summarized by the
following statement.
\begin{theorem}
\label{thm:CAP}
{\bf [Computer-Assisted]}
Consider the value $\varepsilon=0.003$
and the corresponding Hamiltonian
$H^{(0)}$ that has been introduced by
the procedure described in the present
section, starting from the Hamiltonian
$\tilde{H}+f$ which is fully
defined by the
equations~\eqref{eq:hamiltonian1}--\eqref{eq:perturbation0},
\eqref{eq:hamiltonian2}
and~\eqref{eq:control3}.
Then, there exists a torus which is
both invariant with respect to the
Hamiltonian flow induced by $H^{(0)}$
and travelled by quasi-periodic
motions characterized by the
angular velocity vector
$\big(\omega_D\,,\,1\big)$.
\end{theorem}

In order to perform the CAP of the
previous statement, it has been
necessary to execute the algorithm 
discussed in the previous
section~\ref{sec:thm} with
$R_I=150$ and $\RII=8500$. This
required a total computational time
of about 62~hours
on a workstation equipped with 2 CPU
\texttt{Intel XEON-GOLD 5220} (2.2~GHz)
and 384~GB of RAM. Most of the time
has been necessary for the explicit
computation of the (truncated) expansions
of the Hamiltonians $H^{(r)}$ for
$r=1\,,\,\ldots\,,\,R_I=150$; this task
has been made by running the program
\texttt{expl\_transf.c}.

In order to fix the ideas, let us
report here the parameters
entering in the estimates of the Hamiltonian
$H^{(\RII)}$ that has been produced after
the execution of $\RII=8500$ normalization
steps. Their values are
\begin{equation}
\vcenter{\openup1\jot 
\halign{
$\displaystyle{#}\hfil$\qquad&$\displaystyle{#}\hfil$\qquad&$\displaystyle{#}\hfil$\cr
\rho \,=\, 0.024215\,,
&\sigma \,=\, 0.00015686\,,
&\lambda \,=\, 10^{-91.5}\,,
\cr
E \,=\, 0.0061441\,,
&\gamma \,=\, 0.24999\,,
&m \,=\, 1.3809\ .
\cr
}}
\end{equation}
The definitions appearing in
formulae~\eqref{eq:lambda*}--\eqref{eq:paramThm}
allows to straightforwardly compute the
threshold value of $\lambda^*=10^{-85.7}$.
Since the perturbation is so small that
also the crucial hypothesis of
theorem~\ref{thm:kam} is satisfied, i.e.,
$\lambda < \lambda^*$, then the 
program \texttt{iteration+proof.c} claims that
the wanted (invariant) torus exists in
the final part of the output that is
produced at the end of its execution.

\subsection{Determination of $\xi_\text{init}$}
\label{sec:xistar}

Let us consider the Hamiltonian $H_c$
that is defined in equation~\eqref{eq:Hc}.
In particular, we focus on a single
coefficient, i.e.,
\mbox{$B(\xi_{\text{init}})=\langle\frac{\partial H_c}{\partial \psi}\rangle\big|_{\psi=0}\,$}.
We aim to determine the value of the
initial shift on the action $\psi$, i.e.,
$\xi_\text{init}\in\mathbb{R}$ 
so that $B(\xi_\text{init})$ is as near as 
possible to the wanted frequency 
$\omega_\text{init}\,$; for this purpose,
we use a kind of Newton algorithm. We need
a first approximation $\xi_0\,$, that is accurate
enough in order to start a successful
search of $\xi_\text{init}\,$. By comparing 
$e^{\lieDer{\xi_\text{init}\theta}}
e^{\lieDer{X}}\big(\tilde{H}+f\big)$ with the
expansions~\eqref{eq:hamiltonian2}
and~\eqref{eq:control3} of $\tilde{H}$ and $f$,
respectively, one can easily realize that
the main terms that are linear in $\psi$ and
with non-zero angular average are collected in
the following expression:
\begin{equation}
  \tilde{\omega}\psi+
  \,\frac{\partial \langle h_{2}\rangle}{\partial\psi}\,\xi_0\ .
\end{equation}
Therefore, we impose that it is equal to
$\omega_\text{init}\psi+\Pphi$. This allows
us to obtain a linear equation, i.e.,
\begin{equation}
  \tilde{\omega}+\,\frac{\partial^2\langle h_{2}\rangle}{\partial\psi^2}\,\xi_0
  \,=\,\omega_\text{init}\ ,
\end{equation}
that can be easily solved with respect
to the unknown $\xi_0\,$.

The first order expansion of the
equation $B(\xi_{\text{init}})=\omega_\text{init}$
in the linear approximation centered about
$\xi_0$ gives
\begin{equation}
 \omega_\text{init}\,\simeq\,B(\xi_0)
 \,+\,B^{\prime}(\xi_0)\,(\xi_\text{init}-\xi_0)\ .
 \label{eq:implicit_xi}
\end{equation}
This suggest to introduce
a sequence $\{\xi_s\}_{s\geq 0}$ that
is iteratively defined as follows:
\begin{equation}
  \xi_{s+1}\,=\,\xi_s\,+\,\delta\xi_s\,,\quad 
  \delta\xi_s\,=\,\frac{ \omega_\text{init}\, -\, B(\xi_s) }{ \Delta_{\xi_s} B }\ ,
\label{eq:rec_scheme}
\end{equation}
where $\Delta_{\xi_s} B$ is an approximate evaluation
of the derivative $B^{\prime}(\xi_s)$, because
we set
$\Delta_{\xi_s} B=\big[ B\big(\xi_s(1+dx)\big)-B(\xi_s)\big]/(\xi_s\,dx)$,
being $dx$ a positive number such that $dx\ll 1$.
In order to complete the description of the
Newton method, we have also to fix a
tolerance $\Xi$ so that, when 
$|\delta\xi_s|<\Xi$, we stop the algorithm
by setting the value of
$\xi_\text{init}=\xi_{s+1}\,$.
In all our computations dealing with
the construction of the normal form for
the invariant torus characterized by
the angular velocity vector
$\big(\omega_D\,,\,1\big)$,
after having set \mbox{$dx=10^{-3}$}
and \mbox{$\Xi=10^{-8}$},
at most three iterations of the
Newton method have been
enough to determine the value
of $\xi_\text{init}\,$.

Of course, the explicit computations
of the expansions of $e^{\lieDer{\chi_2}}
e^{\lieDer{\xi_s\theta}}
e^{\lieDer{X}}\big(\tilde{H}+f\big)$
are performed by a code using the
facilities for doing algebraic
manipulations that are provided by
the software package
{\it X$\rho$\'o$\nu o\varsigma$}.

\section{Conclusions}
\label{sec:conclusions}

In \cite{chandre_control_2006}, a Hamiltonian
model including a control term was introduced
to study the problem of the
confinement of a strong magnetic field 
with tearing modes. In this work we 
have studied some properties of that 
system by using an approach based on the
Frequency Analysis, and we have also used
these preliminary results as input
for starting a rigorous Computer Assisted
algorithm that is designed to construct
invariant KAM tori.  
Both the Computer Assisted Proof (of which 
we join the whole code) and the preliminary
analysis can be applied with 
small adaptations\footnote{Let us emphasize
that the codes that perform the CAP and are
related to the present work (in the
{\it supplementary material}) can apply
to system with $n>2$ degrees of
freedom, after having modified the
evaluation of the Diophantine constant
$\gamma$. In the attached version of the
codes, that computation is done by using
the continued fraction algorithm. This is
the only constraint restricting the
applications to Hamiltonian models with
$n=2$ degrees of freedom.} to any model
described by a quasi-integrable Hamiltonian
$H(p,q)$ whose Taylor-Fourier expansion
in the action-angle coordinates $(p,q)$ is
finite.

By the Frequency Analysis, we have seen that
the phase space of the perturbed, uncontrolled 
magnetic field is dominated by chaos already 
for $\varepsilon=0.0012$. When we add the 
control term, a regular region appears, and 
it survives for up to $\varepsilon=0.0045$,
nearly 4 times the breaking threshold
of the uncontrolled system. Therefore, the
control term can create a bunch 
of invariant curves, that we may call a 
``transport barrier''; these curves are 
KAM tori. This is in agreement with a known
result obtained by Morbidelli and Giorgilli 
\cite{morbidelli_superexponential_1995}
stating that KAM tori are not isolated,
they fill nearly all the local region in the
vicinity of some ``chief torus'', which has
the highest breaking threshold with respect
to the perturbation parameter $\varepsilon$. 

We have also computed the frequency
related to one of these ``chief
tori'' and we have used our Computer Assisted
Proof to ensure the persistence of the
corresponding invariant surface for
$\varepsilon=0.003$, i.e., up to about
$66$\% of the expected value concerning its
breakdown threshold. Our result is not at
same level of performance with respect to
others existing in the literature. For
instance, more than twenty years ago
in~\cite{celletti_improved_2000}
the existence of a KAM torus was proved
up to values of the small parameter that
are approximately equal to $92$\% of the
numerical breakdown threshold, in the
case of a $1\tfrac{1}{2}$ d.o.f. Hamiltonian
model describing a forced pendulum.
More recently, \cite{figueras_rigorous_2017}
the rigorous proof of the existence of
invariant KAM tori was produced with an
accuracy very close to the optimal value
in the widely studied case of the standard
mapping. However, we think that the problem
considered here is extremely challenging:
as discussed above, due to the presence of
the control term, invariant tori can persist
for much larger values of the perturbation.
Therefore, proving the convergence of the
whole computer assisted procedure is
expected to be much
harder than in problems without any control
term. Applications of approaches based on
other CAP techniques to the Hamiltonian
model considered here could be interesting,
also to compare the performances of
different methods. In this context, let us
emphasize that our approach is in a good
position to obtain rigorous results about
the stability in the vicinity of an
invariant ``chief torus'', because our
algorithm is designed to construct
Kolmogorov normal forms
(see~\cite{giorgilli_kolmogorov_2009}
and~\cite{giorgilli_secular_2017}).

Tearing modes can induce just a minor part
of the zoo of instabilities and transport 
phenomena observed in plasma physics and
concerning the field alone. It would be
interesting to consider system representing
the field interacting with a charged particle. 
The Hamiltonian control theory was applied 
also to a simple model of such a system in 
\cite{ciraolo_control_2004}. An interesting input 
in this direction may come from a new approach
\cite{di_troia_charge_2015} to Guiding Centre 
Theory. It would be even more interesting to 
consider a system of many charged particles in 
mutual interaction. As shown recently 
\cite{carati_transition_2012}, the application
of Hamiltonian techniques to this type of
systems can lead to novel and unexpected
results.

\subsection*{Acknowledgments}
This work was partially supported by the MIUR-PRIN project 20178CJA2B
-- ``New Frontiers of Celestial Mechanics: theory and Applications''
and the ``Beyond Borders'' programme of the University of Rome Tor
Vergata through the project ASTRID (CUP E84I19002250005).  The authors
are indebted with prof. A.~Giorgilli for the availability of the
software package {\it X$\rho$\'o$\nu o\varsigma$} and they acknowledge
also the MIUR Excellence Department Project awarded to the Department
of Mathematics of the University of Rome ``Tor Vergata'' (CUP
E83C18000100006).





\nocite{*}
\bibliography{KAM_plasma_control}
\bibliographystyle{plain}

\end{document}